\documentclass[12pt,twoside]{amsart}
\usepackage{amssymb}
\usepackage{amscd}

\title[Remarks on algebraic fiber spaces] 
{Remarks on algebraic fiber spaces}
\author{Osamu Fujino}
\address{Research Institute for Mathematical Sciences\\ 
Kyoto University, Kyoto 606-8502 Japan}
\email{fujino@kurims.kyoto-u.ac.jp}

\newcommand{\Supp}[0]{{\operatorname{Supp}}}
\newcommand{\xIm}[0]{{\operatorname{Im}}}
\newcommand{\xP}[0]{{\operatorname{Pic^0}}}
\newcommand{\Alb}[0]{{\operatorname{Alb}}}

\newcommand{\xVar}[0]{{\operatorname{Var}}}

\newcommand{\xcd}[0]{{\operatorname{fd_{Alb}}}}
\newcommand{\xd}[0]{{\operatorname{d_{Alb}}}}
\newcommand{\xrank}[0]{{\operatorname{rank}}}
\newcommand{\xG}[0]{{\operatorname{Gal}}}
\newtheorem{thm}{Theorem}[section]
\newtheorem{lem}[thm]{Lemma}
\newtheorem{cor}[thm]{Corollary}
\newtheorem{prop}[thm]{Proposition}

\theoremstyle{definition}
\newtheorem{defn}[thm]{Definition}
\newtheorem{ex}[thm]{Example}
\newtheorem{rem}[thm]{Remark}
\newtheorem*{ack}{Acknowledgements}       

\newtheorem*{notation}{Notation}         
       
\newtheorem{say}[thm]{}
\newtheorem*{claim}{Claim}       

\theoremstyle{remark}

\newtheorem{step}{Step}
\begin{document}
\bibliographystyle{amsalpha+}

\abstract 
In this paper, we collect basic properties of 
the Albanese dimension and explain how to generalize 
the main theorem of \cite{fujino}. 
This paper is a supplement and a generalization of \cite{fujino}. 
We also prove an inequality of irregularities for 
algebraic fiber spaces in the appendix, which is an exposition of 
Fujita-Kawamata's semi-positivity theorem. 
\endabstract

\maketitle
\section{Introduction}

The following is the main theorem of \cite{fujino}. 
In this paper, we recall the Albanese dimension and 
collect its basic properties. 
Then we explain how to generalize Theorem \ref{main} in 
Section \ref{sec4}. 
We note that Section \ref{sec4} is not self-contained. 
Finally, we explain Fujita-Kawamata's semi-positivity 
theorem and prove an inequality of irregularities for 
algebraic fiber spaces in Section \ref{sec3}. 

\begin{thm}[Main theorem of \cite{fujino}]\label{main}
Let $f:X\to Y$ be a surjective morphism 
between non-singular projective varieties with connected fibers. 
Let $F$ be a sufficiently general fiber of $f$. 
Assume that $F$ has maximal Albanese dimension. 
Then $\kappa (X)\geq \kappa (Y)+\kappa 
(F)$. 
\end{thm}

The main theorem of this paper is Theorem \ref{mainda} below. 
See also Theorem \ref{tka}. 
Before we state the main theorem, we 
recall {\em{Albanese dimension}} and 
define {\em{Albanese fiber dimension}} 
for non-singular projective varieties. 
The notion of Albanese dimension was already defined by 
many mathematicians. 

\begin{defn}[Albanese dimension and Albanese fiber dimension]\label{dim}
Let $X$ be a non-singular projective variety. 
Let $\Alb(X)$ be the {\em{Albanese variety}} of $X$ 
and $\alpha_X:X\to \Alb(X)$ the {\em{Albanese map}}. 
We define the {\em{Albanese dimension $\xd(X)$}} as 
follows; 
$$
\xd(X):=\dim (\alpha_X(X)). 
$$ 
We call $\xcd(X):=\dim X-\xd(X)$ the {\em{Albanese fiber 
dimension}} of 
$X$. 
By the definition, it is obvious that $0\leq \xd(X)\leq \dim X$ and 
$0\leq \xcd (X)\leq \dim X$. 
We note that $\xd(X)=\xrank _{\mathcal O_X}(\xIm(H^{0}(X,
\Omega^{1}_{X})\otimes _{\mathbb C}\mathcal O_X\to \Omega^{1}_{X}))$. 

We say that $X$ {\em{has maximal Albanese dimension}} when 
$\dim X=\xd(X)$, equivalently, $\xcd (X)=0$. 
\end{defn}

The following is the main theorem of this paper, 
which is a generalization of Theorem \ref{main}. 
See also Theorem \ref{tka} below. 

\begin{thm}\label{mainda}
Let $f:X\to Y$ be a surjective morphism 
between non-singular projective varieties with connected fibers. 
Let $F$ be a sufficiently general fiber of $f$. 
Assume that $\xcd(F)\leq 3$. 
Then $\kappa (X)\geq \kappa (Y)+\kappa 
(F)$. 
\end{thm}

The following corollary is obvious by Theorem \ref{mainda} and 
Proposition \ref{seisitu} (2) below. 

\begin{cor}\label{4jigen}
Let $f:X\to Y$ be a surjective morphism 
between non-singular projective varieties with connected fibers 
and $F$ a sufficiently general fiber of $f$. 
Assume that $\dim F=4$ and the irregularity $q(F)>0$. 
Then $\kappa (X)\geq \kappa (Y)+\kappa 
(F)$.  
\end{cor}

The idea of the proof of Theorem \ref{mainda} is 
to combine Theorem \ref{main} with Theorem \ref{daiteiri} below. 
It is a special case of \cite[Corollary 1.2]{ka3}. 

\begin{thm}[{c.f.~\cite[Corollary 1.2]{ka3}}]\label{daiteiri}
Let $f:X\to Y$ be a surjective morphism between non-singular 
projective varieties. 
Assume that $\dim X-\dim Y\leq 3$ or sufficiently general 
fiber of $f$ are birationally equivalent to Abelian varieties. 
If $\kappa (Y)\geq 0$, then 
$$ 
\kappa (X)\geq\kappa (F)+\max\{\kappa (Y), \xVar (f)\}, 
$$ 
where $F$ is a sufficiently general fiber of $f$. 
\end{thm}

For $\xVar(f)$ and other positive answers to 
Iitaka's conjecture, see \cite[Sections 6, 7]{M}. 

We summarize the contents of this paper: 
In section \ref{sec2}, we define Albanese dimension and 
Albanese fiber dimension for 
complete (not necessarily non-singular) 
varieties and collect their several basic properties. 
Section \ref{sec4} sketches the proof of the main theorem: 
Theorem \ref{mainda}. 
The proof depends on 
Theorem \ref{main} and the arguments in \cite{fujino}. 
In section \ref{sec3}, which is an appendix, 
we explain Fujita-Kawamata's semi-positivity theorem 
and prove an inequality about {\em{irregularities}} 
for algebraic fiber spaces (Theorem \ref{irr}). 
It may be useful for the study of the relative 
Albanese map. 
The statement is as follows; 

\begin{thm}[Inequality of irregularities]\label{int-irr}
Let $f:X\to Y$ be a surjective morphism between non-singular 
projective varieties with connected fibers. Then 
$$q(Y)\leq q(X)\leq q(Y)+q(F), 
$$ 
where $F$ is a general fiber of $f$. 
\end{thm}

We recommend the readers who are only interested in 
Theorem \ref{int-irr} to read Section \ref{sec3} directly 
after checking the notation below. 
Section \ref{sec3} is independent of the results in Sections 
\ref{sec2} and \ref{sec4}. 
Section \ref{sec3} is intended to the readers who 
are not familiar with the technique of the 
higher dimensional algebraic geometry. 
It is an expository section. 

\begin{ack} 
I talked about parts of this paper at RIMS in June 2002.  
I would like to thank Professors Shigefumi Mori and Noboru 
Nakayama, who pointed out some mistakes. 
I was partially supported by the Inoue Foundation for Science. 
\end{ack}

We fix the notation used in this paper. 

\begin{notation}\label{no} 
We will work over the complex number field $\mathbb C$ 
throughout this paper. 
\begin{itemize}
\item[(i)] 
A {\em{sufficiently general point}} $z$ (resp.~{\em{subvariety}} 
$\Gamma$) of the variety $Z$ 
means that $z$ (resp.~$\Gamma$) is not contained in the countable union of 
certain proper Zariski closed subsets. 
We say that a subvariety $\Gamma$ (resp.~point $z$) 
is {\em{general}} in $Z$ if 
it is not contained in a 
certain proper Zariski closed subsets. 

Let $f:X\to Y$ be a morphism between varieties. 
A {\em{sufficiently general fiber}} (resp.~{\em{general fiber}}) 
$X_y=f^{-1}(y)$ of $f$ means 
that $y$ is a sufficiently general point (resp.~general) in $Y$. 
\item[(ii)] 

An {\em{algebraic fiber space}} $f:X\to Y$ is a proper surjective 
morphism between non-singular projective varieties 
$X$ and $Y$ with connected fibers. 

\item[(iii)] 
Let $f:X\to Y$ be a surjective morphism between varieties. 
We put $\dim f:=\dim X-\dim Y$. 
Let $X$ be a variety and $\mathcal F$ a
coherent sheaf on $X$. 
We write $h^i(X,\mathcal F)=\dim _{\mathbb C} H^i(X,\mathcal F)$. 
If $X$ is non-singular projective, then $q(X):=h^1(X,\mathcal O_X)$ 
denotes the {\em{irregularity}} of $X$. 

\item[(iv)] The words {\em{locally free sheaf}} and 
{\em{vector bundle}} are used interchangeably. 

\item[(v)] Since most questions we are interested 
in are birational ones, we usually make birational 
modifications freely whenever it is necessary. 
If no confusion is likely, we denote the new objects with 
the old symbols. 

\item[(vi)]
Let $X$ be a non-singular projective variety. 
If the {\em{Kodaira dimension}} $\kappa(X)>0$,
then we have the {\em{Iitaka fibration}} $f: X\to Y$, 
where $X$ and $Y$ are non-singular projective varieties 
and $Y$ is of dimension $\kappa(X)$, such that the sufficiently 
general fiber of $f$ is non-singular,
irreducible with $\kappa =0$.  
Since the Iitaka fibration is determined only up to
birational equivalence, we used the above abuses in (v). 
For the basic properties of the Kodaira dimension and 
the Iitaka fibration, 
see \cite[Chapter III]{ueno} or \cite[Sections 1,2]{M}. 
\end{itemize}
\end{notation}

\section{Albanese dimension}\label{sec2}

In this section, we collect several basic properties of the 
Albanese dimension and the Albanese fiber dimension. 
The next lemma is easy to check. 

\begin{lem}\label{birinv}
Let $f:W \to V$ be a birational morphism 
between non-singular projective varieties. 
Then $(\Alb(V), \alpha_V\circ f)$ is the 
Albanese variety of $W$. 
In particular, 
$\xd(V)=\xd(W)$.   
\end{lem}
\proof
This is obvious. 
See, for example, \cite[Proposition 9.12]{ueno}.  
\endproof

By the above lemma, we can define the Albanese dimension and 
the Albanese fiber dimension for singular varieties. 

\begin{defn}\label{honto}
Let $X$ be a complete variety. 
We define $\xd(X):=\xd(\widetilde X)$, 
where $\widetilde X$ is a non-singular projective 
variety that is birationally equivalent to 
$X$. 
We put $\xcd(X):=\dim X-\xd(X)$. 

By Lemma \ref{birinv}, the {\em{Albanese dimension $\xd(X)$}} 
and the {\em{Albanese fiber dimension 
$\xcd(X)$}} are well-defined 
birational invariants of $X$. 
\end{defn}

\begin{defn}[Varieties of maximal Albanese dimension]\label{mad}
Let $X$ be a complete variety. 
If $\xcd (X)=0$, 
then we say that $X$ is {\em{of maximal 
Albanese dimension}} or {\em{has maximal 
Albanese dimension}}. 
We note that $\kappa (X)\geq 0$ if $\xcd (X)=0$. 
By \cite[Theorem 1]{ka0}, 
$\kappa (X)=0$ and $\xcd (X)=0$ 
implies that $\alpha _X:X\to \Alb (X)$ is birational. 
\end{defn}

\begin{defn}[Albanese dimension of fibers]\label{fib}
Let $f:X\to Y$ be a surjective morphism between 
non-singular projective varieties. 
By taking the Stein factorization of $f$, we obtain; 
$$
\begin{matrix}
X\ \ & \longrightarrow& \ \ Z \\
{\ \ \searrow} & \ &  {\swarrow\ \ } \\
 \ & Y & \ & .
\end{matrix}
$$ 
We shrink $Z$ suitably and take a finite \'etale cover 
$\widetilde Z\to Z$ such that $\widetilde X:=X\times _Z 
\widetilde Z\to \widetilde Z$ has a section. 
Then, we obtain a relative Albanese map $\widetilde X\to 
\Alb(\widetilde X/\widetilde Z)\to \widetilde Z$ (see also 
\ref{relalb} below). 
By this relative Albanese map, it is easy to see that 
the Albanese dimension $\xd(X_z)$ is 
independent of $z\in U$, 
where $U$ is a suitable non-empty Zariski open 
set of $Z$ and $X_z$ is a fiber of $X\to Z$. 
We note that we can assume that 
the irregularity $q(X_z)$ is independent 
of $z\in U$. 

We put $\xd(f):=\xd(X_z)$, where $z$ is a general point of $Z$ and 
$\xcd (f):=\dim f-\xd(f)$. 
We call $\xd(f)$ (resp.~$\xcd(f)$) the {\em{Albanese dimension}} 
(resp.~{\em{Albanese fiber dimension}}) of $f$. 
\end{defn}

The following proposition will play important roles 
in the proof of the main theorem 
(c.f.~\cite[Proposition 2.3 (3)]{fujino}).  

\begin{prop}\label{sub}
Let $V$ be a non-singular projective variety 
and $X$ a closed subvariety of $V$. 
If $X$ is general in $V$, then 
$$
\xcd(X)\leq \dim X-\dim \alpha_V(X)\leq \xcd (V). 
$$ 
\end{prop}
\proof 
Let $\widetilde X\to X$ be a resolution. 
We consider the following commutative diagram; 
$$
\begin{CD}
\widetilde X &@>>> & \Alb (\widetilde X)\\
@VVV& & @VVV \\
V &@>{\alpha _V}>> & \Alb (V).  
\end{CD}
$$ 
Since $X$ is general in $V$, $\dim V-\dim X\geq 
\xd(V)-\dim \alpha_V(X)$. 
Therefore, $\xcd(V)\geq \dim X-\dim \alpha_V(X)\geq \xcd (X)$. 
\endproof

The next example shows that the equality does not 
necessarily holds in Proposition \ref{sub}. 

\begin{ex}\label{kantan}
Let $X:=\mathbb P^1\times E$, 
where $E$ is an elliptic curve. 
Let $p:X\to E$ be the second projection and 
$C$ an irreducible curve on $X$. 
Then $\xcd (C)=0$ unless $p(C)$ is a point. 
If $C$ is a fiber of $p$, 
then it is obvious that $\xcd (C)=1$. 
We note that $\xcd (X)=1$. 
\end{ex}

\begin{ex}
Let $A$ be an Abelian surface and $X$ a one point blow-up 
of $A$. Let $E$ be the $(-1)$-curve on $X$. 
Then $\xcd (E)=1$ since $E\simeq \mathbb P^1$. 
Therefore, $1=\xcd (E)>\xcd (X)=\xcd (A)=0$.  
\end{ex}

\begin{cor}[Easy addition of the Albanese dimension]\label{eas}
Let $f:X\to Y$ be a surjective morphism between non-singular 
projective 
varieties. 
Then 
$$\xd(X)\leq \dim \alpha _X(F)+\dim Y \leq \xd(f)+\dim Y. 
$$ 
\end{cor}
\proof
This easily follows from Proposition \ref{sub}. 
We note that Definition \ref{fib}. 
\endproof

\begin{ex}\label{new}
Let $A$ be an Abelian surface and $H$ be a non-singular very ample 
divisor on $A$. 
We take a general member $H'\in |H|$. 
We can write $H'=H+(h)$, where $h\in \mathbb C(A)$. 
Consider the rational map; 
$$
h:A\dashrightarrow \mathbb P^1. 
$$ 
By blowing up the points $H\cap H'$, 
we obtain an algebraic fiber space; 
$$
f:X\longrightarrow Y, 
$$ 
which is birationally equivalent to $h$. 
In this case, $\xd(X)=2$, $\xd(Y)=0$, and $\xd(f)=1$. 
Therefore, we can not replace $\dim Y$ with $\xd(Y)$ in Corollary 
\ref{eas}. 
\end{ex}

The following claim is a variant of \cite[Proposition 2.4]
{fujino}.  

\begin{prop}\label{iifiber}
Let $V$ be a non-singular projective variety and 
$f:V\to W$ be the Iitaka fibration. 
Then $\xcd (W)\leq \xcd (V)$. 
In particular, if $V$ is of maximal Albanese dimension, 
then so is $W$. 
\end{prop}
\proof
We put $m:=\xcd (V)$ and $n:=\dim V$. 
Let $F$ be a sufficiently general fiber of $f$. 
Then $\kappa (F)=0$ and 
$\xcd(F)\leq m$ by Proposition \ref{sub}. 
By \cite[Theorem 1]{ka0}, 
$\alpha_F:F\to \Alb(F)$ is an algebraic fiber space. 
By the following diagram; 
$$
\begin{CD}
F &@>>> & \Alb (F)\\
@VVV& & @VVV \\
V &@>{\alpha _V}>> & \Alb (V), 
\end{CD}
$$ 
$\alpha_V(F)$ is an Abelian variety. 
Since there exists at most countably many Abelian subvarieties 
in $\Alb(V)$, 
there is an Abelian subvariety $A$ of $\Alb(V)$ 
such that $\alpha_V(F)$ is a translation of $A$ 
for general fibers $F$. 
We note that $\dim A\leq \dim F=n-\dim W=n-\kappa (V)$. 

Let $\psi: \Alb(V)\to \Alb(V)/A$ be the quotient map. 
By the definition of $A$, 
$\psi\circ \alpha_V$ induces a rational map 
$\varphi:W\dasharrow \Alb(V)/A$. 
Since $\Alb(V)/A$ is Abelian, $\varphi$ is a morphism. 
By the universality of $(\Alb(W), \alpha_W)$, 
$\varphi$ factors through $\Alb(W)$. 
Therefore, 
\begin{eqnarray*}
\xd(W)&\geq&\dim \varphi(W)\\
&=&\dim \psi(\alpha_V(V))\\
&\geq& n-m-\dim F\\
&=& \kappa(V)-m.
\end{eqnarray*} 
Note that $\dim W=\kappa (V)$. 
Thus, we have the required inequality $\xcd (W)\leq m$. 
\endproof

\begin{rem}\label{iifiber2}
We can generalize Proposition \ref{iifiber} as follows without 
difficulties. Details are left to the readers. 

Let $V\to W$ be a surjective morphism between non-singular projective 
varieties. 
If $q(f)=\xd(f)$, then $\xcd (W)\leq\xcd(V)$. 
\end{rem}

The following example says that the equality doesn't necessarily 
hold in Proposition \ref{iifiber}. 

\begin{ex}\label{ex2}
Let $S$ be a $K3$ surface and $C$ is a non-singular projective 
curve with the genus $g(C)\geq 2$. 
We put $X:=S\times C$. 
Then the second projection $X\to C$ is the Iitaka 
fibration. 
It is easy to check that $\xcd (X)=2$ and $\xcd (C)=0$.  
\end{ex}

We collect several basic properties of the Albanese dimension 
and Albanese fiber dimension 
for the reader's convenience. First, we recall the 
following obvious fact. 

\begin{lem}\label{akiraka}
Let $X$ be a non-singular projective variety. 
Then the irregularity $q(X)=0$ if and only if $\xd(X)=0$. 
\end{lem}

\begin{prop}\label{seisitu}
{\em{(1)}} Let $X$ be a complete variety with 
$\dim X=n$. 
Then $\xcd (X)\leq n$. 

{\em{(2)}} Let $X$ be an $(n+1)$-dimensional 
non-singular projective variety with the irregularity 
$q(X)>0$. 
Then $\xcd (X)\leq n$. 

{\em{(3)}} Let $f:X\to Y$ be a generically finite morphism between 
non-singular projective varieties. 
Then $\xd (X)\geq \xd (Y)$. 
Equivalently, $\xcd (X)\leq \xcd (Y)$.  

{\em{(4)}} Let $f:X\to Y$ be a surjective morphism between 
non-singular projective varieties. 
Then $\xcd (X)\leq \xcd(Y)+\dim f$. 

{\em{(5)}} Let $X$ and $Y$ be non-singular projective varieties. 
Then $\xd(X\times Y)=\xd(X)+\xd(Y)$. 
Equivalently, $\xcd(X\times Y)=\xcd(X)+\xcd(Y)$. 

\proof 
The claims (1), (3), and (4) are obvious by the definition. 
For (2), we note that $\alpha _X(X)$ is not a point by Lemma 
\ref{akiraka}. 
For (5), note that $(\Alb(X)\times \Alb(Y), \alpha_X\times \alpha_Y)$ 
is the Albanese variety of $X\times Y$ by 
the K\"unneth formula. 
\endproof 
\end{prop}

\begin{lem}\label{lef}
Let $V$ be a non-singular projective variety 
and $H$ a non-singular ample divisor on $V$. 
Then $\Alb(H)\to \Alb(V)$ is 
surjective {\em{(}}resp.~an isomorphism{\em{)}} 
if $\dim V\geq 2$ 
{\em{(}}resp.~$\dim V\geq 3${\em{)}}. 
\end{lem}
\proof
Apply the Kodaira vanishing theorem to 
$H^1(V,\mathcal O_V(-H))$ and 
$H^2(V,\mathcal O_V(-H))$. 
\endproof

\begin{prop}\label{lef2}
Let $V$ be a non-singular projective variety 
and $H$ a non-singular ample divisor on $V$. 
Assume that $H$ is general in $V$. 
Then $\xcd(H)=\xcd (V)-1$ if $\dim V\geq 3$. 
\end{prop}
\proof
It is obvious by Lemma \ref{lef}. 
\endproof

The following lemma is a key lemma in Section \ref{sec4} 
(see \ref{relalb} below). 

\begin{lem}\label{fiber}
Let $X$ be a non-singular projective variety and $\alpha_X: 
X\to \Alb(X)$ the Albanese mapping. 
We consider the Stein factorization; 
$$
X\longrightarrow W\longrightarrow\alpha_X(X). 
$$ 
Then $\xcd (W)=0$, that is, $W$ has maximal Albanese dimension. 
\end{lem}
\proof
It is obvious by the definition. 
See Proposition \ref{seisitu} (3). 
\endproof

\section{Sketch of the proof of the main theorem}\label{sec4}

In this section, we sketch the proof of the main theorem: 
Theorem \ref{mainda}. 
We explain how to modify the arguments in \cite[Section 4]{fujino}. 
For the details, see \cite{fujino}. 

\begin{say}[Algebraic fiber space associated to the relative Albanese 
map]
\label{relalb}
Let $f:X\to Y$ be an algebraic fiber space. 
From now on, we often replace $Y$ with its non-empty Zariski open 
set and $X$ with its inverse image. 
We denote the new objects with the old symbols. 
We can assume that $f$ is smooth. 
We can assume that there exists a finite 
\'etale cover $\pi:\widetilde Y\to Y$ such 
that $\pi$ is Galois and $\widetilde f:\widetilde X:=
X\times _Y\widetilde Y\to \widetilde Y$ has a 
section. Then there exists the relative Albanese map 
$\alpha _{\widetilde X/\widetilde Y}:\widetilde X\to \Alb(\widetilde X/
\widetilde Y)$ over $\widetilde Y$. 
We note that $\Alb(\widetilde X/
\widetilde Y)\simeq \xP(\xP (\widetilde X/\widetilde Y)/\widetilde Y)$. 
Put $k:=\deg\pi$ and $G:=\xG(\widetilde Y/Y)$ the 
Galois group of $\pi$. 
Then $G$ acts on $\widetilde X=X\times _Y\widetilde Y$ and 
$\Alb (\widetilde X/\widetilde Y)=\Alb(X/Y)\times _Y\widetilde Y$. 
Thus $G$ acts on $\alpha _{\widetilde X/\widetilde Y}$ 
as follows; 
$$
(g\cdot \alpha _{\widetilde X/\widetilde Y})(x)=
g^{-1}(\alpha _{\widetilde X/\widetilde Y}(gx)), 
$$ 
where $g\in G$ and $x\in \widetilde X$. 
We put $\beta_{\widetilde X/\widetilde Y}:=
\sum _{g\in G}g\cdot \alpha _{\widetilde X/\widetilde Y}: 
\widetilde X\to \Alb(\widetilde X/\widetilde Y)$. 
Let $\widetilde \Gamma \subset \widetilde X\times _{\widetilde Y}
\Alb(\widetilde X/\widetilde Y)$ be the graph of 
$\beta_{\widetilde X/\widetilde Y}$ and 
$\Gamma \subset X\times _{Y}
\Alb(X/Y)$ the image of $\widetilde \Gamma$. 
By the construction of $\beta_{\widetilde X/\widetilde Y}$, 
$\Gamma$ induces a morphism 
$\beta _{X/Y}:X\to \Alb(X/Y)$ over $Y$. 
Let us see $\beta_{X/Y}$ fiberwise. 
Then it is the composition of the Albanese map and 
the multiplication by $k$ 
of the Albanese variety up to a translation. 
By taking the Stein factorization of $\beta_{X/Y}$, 
we obtain $X\to Z\to Y$. Compactify $X$, $Y$, and $Z$. 
Then, after taking resolutions, 
we obtain; 
$$
\begin{CD}
f:X@>{g}>>Z@>{h}>>Y       
\end{CD}
$$ 
such that 
\begin{itemize}
\item[(i)] it is birationally equivalent to the 
given fiber space $f:X\to Y$. 
\item[(ii)] $X$, $Y$, and $Z$ are non-singular projective 
varieties. 
\item[(iii)] $g$ and $h$ are algebraic fiber spaces. 
\item[(iv)] the general fibers of $h$ have maximal Albanese 
dimension. 
\item[(v)] if $\kappa (X_y)\geq 0$ for sufficiently 
general fibers $X_y$ of $f$, then $\kappa (X_z)\geq 0$ 
for sufficiently general points $z\in Z$. 
This is an easy consequence of the easy addition of the Kodaira dimension. 
\end{itemize}
\end{say}

\proof[Proof of {\em{Theorem \ref{mainda}}}] 
Let $f:X\to Y$ be the given fiber space. 

\begin{step}
If $\kappa (Y)=-\infty$ or $\kappa (F)=-\infty$, 
then the inequality is obviously true. 
From now on, we assume that $\kappa (Y)\geq 0$ and 
$\kappa (F)\geq 0$. 
\end{step}

\begin{step}\label{s2}
If $\kappa (F)=0$, then we construct the fiber space associated to 
the relative 
Albanese map \ref{relalb}, 
we obtain 
$$
\begin{CD}
f:X@>{g}>>Z@>{h}>>Y,   
\end{CD}
$$ 
such that $\dim X-\dim Z\leq 3$ by 
the assumption $\xcd(F)\leq 3$ and the sufficiently general fiber of 
$h$ are of maximal Albanese dimension by \ref{relalb} (iv). 
By \ref{relalb} (v), the Kodaira dimension 
of the sufficiently general fiber of $g$ is non-negative. 
Therefore, by Theorem \ref{main} and Theorem \ref{daiteiri}, 
\begin{eqnarray*}
\kappa (X)&\geq &\kappa (Z)+\kappa (X_z)\\
&\geq&\kappa (Z)\\
&\geq&\kappa (Y)+\kappa (Z_y)\\
&\geq& \kappa(Y), 
\end{eqnarray*} 
where $z$ (resp.~$y$) is a sufficiently general point of $Z$ 
(resp.~$Y$). 
Note that $\kappa (Z_y)\geq 0$ since $\xcd (Z_y)=0$ 
(see Definition \ref{mad}). 
Therefore, we obtained the required inequality when 
$\kappa (F)=0$. 
From now on, we assume that $\kappa (F)>0$. 
\end{step}

\begin{step} 
We recall the following useful lemma (\cite[Lemma 4.2]{fujino}). 

\begin{lem}[Induction Lemma]\label{induction}
Under the same notation as in {\em{Theorem \ref{mainda}}}, 
it is sufficient to prove that 
$\kappa (X)>0$ on the assumption 
that 
$\kappa (Y)\geq 0$ and $\kappa (F)>0$. 
\end{lem}
\begin{rem}\label{tyuui}
We prove this lemma in \cite{fujino} on the assumption that 
$\xcd (F)=0$. The same proof works on the 
weaker assumption that $\xcd (F)\leq 3$. 
We only have to replace the words ``maximal Albanese dimension'' 
with ``$\xcd(\cdot)\leq 3$'' in the proof of \cite[Lemma 4.2]{fujino}. 
The key point is Proposition \ref{sub}. By this, the inductive 
arguments work. 
Details are left to the reader. 
\end{rem}
\end{step}

\begin{step}
By taking the fiber space associated to the relative Albanese 
map \ref{relalb}, we obtain 
$$
\begin{CD}
f:X@>{g}>>Z@>{h}>>Y.    
\end{CD}
$$ 
Since we want to prove $\kappa (X)>0$, we can assume that 
the sufficiently general fibers of $g$ and $h$ have 
zero Kodaira dimension by Theorem \ref{main} and Theorem 
\ref{daiteiri} (see Step \ref{s2} above). 
On this assumption, the general fibers of 
$h$ are birationally equivalent to Abelian varieties. 
We can further assume that $\xVar(g)=\xVar(h)=0$ by 
Theorem \ref{daiteiri}. 
Therefore, we can apply the same proof as in \cite{fujino}. 
Then we obtain $\kappa (X)>0$. 
For the details, see the latter part of \cite[Proof of the theorem]{fujino}. 
We note that $\kappa (F)\geq 1$.  
\end{step}
Therefore, we complete the proof. 
\endproof

The following theorem is obvious by the proof of Theorem \ref{mainda}. 
For the conjecture $C^+_{n,m}$, see \cite[Section 7]{M} and 
Theorem \ref{daiteiri}. 

\begin{thm}\label{tka}
Suppose that $C^+_{n,m}$ holds for every algebraic fiber spaces 
on the assumption that $n-m\leq k$. 
Let $f:X\to Y$ be an algebraic fiber space. 
If $\xcd(F)\leq k$, where $F$ is a sufficiently general 
fiber of $f$, then $\kappa (X)\geq \kappa (Y)+\kappa (F)$. 
\end{thm}

\section{Appendix: Semi-positivity and Inequality of 
Irregularities}\label{sec3}

The aim of this section is to explain Fujita-Kawamata's 
semi-positivity 
theorem in the geometric situation, 
that is, we prove the semi-positivity 
of $R^if_*\omega_{X/Y}$ on suitable assumptions, 
where $f:X\to Y$ is a surjective morphism between non-singular 
projective varieties. 
Then we prove an inequality of irregularities 
for algebraic fiber spaces (see Theorem \ref{irr}), 
which is an easy consequence of the semi-positivity 
theorem. 
This inequality seems to be useful when we 
treat (relative) Albanese variety. 

In spite of its importance, the results and the 
statements about the semi-positivity theorem 
are scattered over various papers (see \cite[\S 5]{M}). 
This is one of the reason why I decided to write down this section
\footnote{It took long time to find the 
statement \cite[Theorem 2]{ka2}. 
I hope that this section will 
contribute to distribute Fujita-Kawamata's semi-positivity 
package.}. 
It is surprising that there are no good references 
about the semi-positivity of $R^if_*\omega_{X/Y}$. 
Kawamata's proof of the semi-positivity theorem 
(c.f.~\cite[Theorem 2]{ka2} and \cite[\S 4]{ka0}) 
heavily relies on the asymptotic behavior of the 
Hodge metric near a puncture. 
It is not so easy for the non-expert to 
take it out from \cite[\S 6]{sch}. 
We recommend the reader to see \cite[Sections 2, 3]{pet} or 
\cite{zu}\footnote
{For the proof of the semi-positivity theorem \cite[\S 4 (2)]{ka2}, 
it is sufficient to know the asymptotic behavior of 
the Hodge metric of the VHS on a {\em{curve}}.}. 
Our proof depends on 
\cite[Proposition 1]{ka2}, \cite{koI}, \cite{koll}, and 
Viehweg's technique. It is essentially the 
same as \cite[Corollary 3.7]{koI}. 
It is much simpler than the original proof. 
We note that Kawamata's proof can be applied to non-geometric 
situations. So, his theorem is much stronger than the results 
explained in this section. See the original article \cite[\S 4]{ka0}. 

\begin{rem} 
This section is a supplement of \cite[\S 5 Part I]{M}, 
especially, \cite[(5.3) Theorem]{M}. 
The joint paper with S.~Mori \cite{FM} generalized 
\cite[\S 5 Part II]{M} and treated several applications. 
The paper \cite{fuji} gave a precise proof of 
\cite[(5.15.9)(ii)]{M} from the Hodge theoretic viewpoint.  
\end{rem}

Let us recall the definition of {\em{semi-positive vector bundles}}. 

\begin{defn}[Semi-positive vector bundles]\label{sp}
Let $V$ be a complete variety and $\mathcal E$ a 
locally free sheaf on $V$. 
We say that $\mathcal E$ is {\em{semi-positive}} 
if and only if the tautological line bundle 
$\mathcal O_{\mathbb P_V(\mathcal E)}(1)$ 
is nef on $\mathbb P_V (\mathcal E)$. 
We note that $\mathcal E$ is semi-positive if and 
only if for every complete curve $C$ and 
morphism $g:C\to V$ every quotient line bundle of $g^*\mathcal E$ 
has non-negative degree. 
\end{defn}

The following result is obtained by 
Koll\'ar and Nakayama. 
For the details, see, for example, the original articles 
\cite[Theorem 2.6]{koll}, \cite[Theorem 1]
{nakayama}, or \cite[(5.3), (5.4)]{M}. 

\begin{thm}\label{hs}
Let $f:V\to W$ be a projective surjective morphism 
between non-singular varieties. 
We assume that there exists a simple normal crossing 
divisor $\Sigma$ on $W$ such that $f$ is smooth 
over $W\setminus \Sigma$. 
Then $R^i f_*\mathcal O_V$ and $R^if_*\omega_{V/W}$ are 
locally free for every $i$. 
Put $W_0:=W\setminus \Sigma$, $V_0:=f^{-1}(W_0)$, $f_0:=f|_{V_0}$, 
and $k:=\dim f$. 
Moreover, if all the local monodromies on 
$R^{k+i}f_{0*}\mathbb C_{V_0}$ around $\Sigma$ are unipotent, 
then $R^if_*\omega_{V/W}$ is characterized by the {\em{canonical 
extension}} of $R^{k+i}f_{0*}\mathbb C_{V_0}$. 
\end{thm}

The following is the main theorem of this section (see \cite[(5.3) 
Theorem]{M}). 

\begin{thm}[Semi-positivity theorem]\label{semi}
Let $f:V\to W$ be a surjective morphism between 
non-singular projective varieties 
with $k:=\dim f$. 
Let $\Sigma$ be a simple normal crossing divisor on $W$ such that 
$f$ is smooth over $W_0:=W\setminus \Sigma$. Put 
$V_0:=f^{-1}(W_0)$ and $f_0:=f|_{V_0}$. 
We assume that all the local monodromies on $R^{k+i}f_{0*}\mathbb 
C_{V_0}$ around $\Sigma$ are unipotent. 
Then $R^if_*\omega_{V/W}$ is a semi-positive vector bundle on 
$W$. 
\end{thm}

Before we prove Theorem \ref{semi}, 
we fix the notation and convention used below. 

\begin{say}
Let $f:V\to W$ be a surjective morphism between 
varieties. 
Let $$
f^s:V^s:=V\times _W V\times_W\cdots \times _W V \to W 
\ \ {\text{(product taken $s$ times)}}
$$
$$
V^{(s)}={\text{desingularization of}} \ \ 
V^s, f^{(s)}:V^{(s)}\to W. 
$$ 
\end{say}

\begin{lem}\label{vie}
On the same assumption as in {\em{Theorem \ref{semi}}}, we have that 
$(R^if_*\omega_{V/W})^{\otimes s}$ is a direct summand 
of $R^{si}f^{(s)}_*\omega_{V^{(s)}/W}$ for 
every positive integer $s$. 
We note that $(f_*\omega_{V/W})^{\otimes s}\simeq 
f^{(s)}_*\omega_{V^{(s)}/W}$ for every positive integer $s$. 
\end{lem}
\proof
We use the induction on $s$. 
First, when $s=1$, the claim is obvious. 
Next, we assume that $(R^if_*\omega_{V/W})^{\otimes (s-1)}$ 
is a direct summand 
of $R^{(s-1)i}f^{(s-1)}_*\omega_{V^{(s-1)}/W}$. 
We consider the following commutative diagram; 
$$
\begin{CD}
V@<<<V\times _W V^{(s-1)} @<<<V^{(s)}\\ 
@V{f}VV @VVV @VV{g}V \\
W@<<{f^{(s-1)}}< V^{(s-1)}@<<{id}<V^{(s-1)}. 
\end{CD}
$$ 
We can assume that $(f^{(s-1)})^{-1}(\Sigma)$ is simple normal 
crossing without loss of generality. 
Since $R^if_*\omega_{V/W}$ is the canonical extension of 
$R^if_{0*}\omega_{V_0/W_0}$ by Theorem \ref{hs} 
(see \cite[Theorem 2.6]{koll}, 
\cite[Theorem 1]{nakayama}), we have the following isomorphism: 
$f^{(s-1)*}R^if_*\omega_{V/W}\simeq R^ig_*\omega _{V^{(s)}/V^{(s-1)}}$ 
by \cite[Proposition]{ka2}. 
Therefore, 
\begin{eqnarray*}
R^jf^{(s-1)}_*R^ig_*\omega _{V^{(s)}/W}&\simeq& 
R^jf^{(s-1)}_*R^ig_*(\omega_{V^{(s)}/V^{(s-1)}}\otimes 
g^*\omega _{V^{(s-1)}/W})\\ 
&\simeq&
R^jf^{(s-1)}_*(R^ig_*\omega_{V^{(s)}/V^{(s-1)}}\otimes 
\omega _{V^{(s-1)}/W}) \\
&\simeq&
R^jf^{(s-1)}_*(f^{(s-1)*}R^if_*\omega_{V/W}\otimes 
\omega _{V^{(s-1)}/W})\\ 
&\simeq&
R^if_*\omega _{V/W}\otimes R^jf^{(s-1)}_*\omega_{V^{(s-1)}/W},  
\end{eqnarray*} 
for every $j$. 
By \cite[Theorem 3.4]{koll}, $R^{(s-1)i}f^{(s-1)}_*R^ig_* \omega 
_{V^{(s)}/W}$ is a direct summand of $R^{si}f^{(s)}_*\omega 
_{V^{(s)}/W}$. 
We assumed that $(R^if_*\omega_{V/W})^{\otimes (s-1)}$ 
is a direct summand 
of $R^{(s-1)i}f^{(s-1)}_*\omega_{V^{(s-1)}/W}$. 
Therefore, $(R^if_*\omega_{V/W})^{\otimes s}$ 
is a direct summand 
of $R^{si}f^{(s)}_*\omega_{V^{(s)}/W}$. 
\endproof

The following lemma is obvious by Koll\'ar's vanishing theorem 
\cite[Theorem 2.1]{koI}. 
For the {\em{regularities}}, see \cite[p.307 Definitions 1, 2 and 
Proposition 1]{klei}. 

\begin{lem}\label{reg}
Let $f:X\to Y$ be a surjective morphism between projective varieties, 
$X$ non-singular, $m:=\dim Y+1$. Let $L$ be an ample line 
bundle on $Y$ which is generated by its global sections. 
Then $R^if_*\omega_X$ is $m$-regular 
{\em{(}}with respect to $L${\em{)}} for every $i$. 
That is, $H^j(Y,R^if_*\omega_X\otimes L^{\otimes (m-j)})=0$ for 
every $j>0$. 
In particular, $R^if_*\omega_X\otimes L^{\otimes l}$ is 
generated by its global sections for $l\geq m$.  
\end{lem}

The next lemma is not difficult to prove. 
For the proof, see, for example, \cite[Proof of Corollary 3.7]{koI}. 

\begin{lem}\label{kolem}
Let $\mathcal E$ be a vector bundle on a complete 
variety $Y$. Assume that there exists a line bundle $\mathcal L$ 
on $Y$ and $\mathcal L\otimes \mathcal E^{\otimes s}$ 
is generated by its global sections for every $s$. 
Then $\mathcal E$ is semi-positive. 
\end{lem}
\proof[Sketch of the proof] 
It is not difficult to see that 
$\mathcal O_{\mathbb P_Y(\mathcal E)}(s)\otimes \pi ^*\mathcal L$ 
is nef for every $s$, where $\pi:\mathbb P_Y(\mathcal E)\to Y$. 
Therefore, $\mathcal O_{\mathbb P_Y(\mathcal E)}(1)$ is 
nef. 
\endproof

\proof[Proof of {\em{Theorem \ref{semi}}}] 
Let $L$ be an ample line bundle on $W$ which is 
generated by its global sections. 
We put $\mathcal E:=R^if_*\omega_{V/W}$ and 
$\mathcal L:=L^{\otimes m}\otimes \omega_W$, where $m=\dim Y+1$. 
Then $R^{si}f^{(s)}_*\omega_{V^{(s)}/W}\otimes \mathcal L$ 
is generated by its global sections for every $s$ by Lemma \ref{reg}. 
By Lemma \ref{vie}, 
$(R^{i}f_*\omega_{V/W})^{\otimes s}\otimes \mathcal L$ 
is generated by its global sections. Therefore, 
$R^{i}f_*\omega_{V/W}$ is semi-positive by Lemma \ref{kolem}. 
\endproof

The following theorem is well-known. 
We write it for the reader's convenience. 
We will use it in the proof of Theorem \ref{irr}. 
We reduce it to Theorem \ref{semi} 
by the semi-stable reduction theorem. 

\begin{thm}\label{curve}
In {\em{Theorem \ref{semi}}}, if $W$ is a curve, then 
the assumption on the monodoromies is not necessary, 
that is, $R^if_*\omega_{V/W}$ is always semi-positive 
for every $i$. 
\end{thm}
\proof 
Without loss of generality, we can assume that $\Supp f^*(P)$ 
is simple normal crossing for every point $P\in W$. 
By the semi-stable reduction theorem, 
we consider the following commutative diagram;  
$$
\begin{CD}
V@<<<V\times _W\widetilde W@<{\nu}<<V'@<<<\widetilde V\\ 
@V{f}VV @VV{g}V @VV{f'}V @VV{\widetilde f}V\\ 
W@<<{\pi}<\widetilde W@<<{id}<\widetilde W@<<{id}<\widetilde W, 
\end{CD}
$$ 
where $\pi:\widetilde W\to W$ is a finite cover, 
$\nu$ is the normalization, and 
$\widetilde f:\widetilde V\to\widetilde W$ is 
a semi-stable reduction of $f$. 
By the flat base change theorem, 
we have $\pi^*R^if_*\omega_{V/W}\simeq 
R^ig_*\omega_{V\times _W\widetilde W/\widetilde W}$ for 
every $i$. 
We consider the following exact sequence, 
which is induced by the trace map;  
$$
\begin{CD}
0@>>> \nu_*\omega_{V'/\widetilde W}@>{tr}>> 
\omega_{V\times_W \widetilde W/\widetilde W}@>>>\delta@>>>0, 
\end{CD}
$$ 
where $\delta$ is the cokernel of $tr$. 
Since $\Supp\, R^{i-1}g_*\delta\subsetneq \widetilde W$, 
we obtain a generically isomorphic inclusion; 
$$
0\to R^if'_*\omega_{V'/\widetilde W}
\simeq R^i\widetilde f_*\omega_{\widetilde V/\widetilde W}\to 
R^ig_*\omega_{V\times _W\widetilde W/ \widetilde W}
\simeq 
\pi^*R^if_*\omega _{V/W}. 
$$ 
We note that $\nu$ is finite and 
$V'$ has at worst rational Gorenstein singularities. 
Since $\widetilde f$ is semi-stable, all the local monodromies 
are unipotent (see, for example, \cite[(4.6.1)]{M}). 
Thus $R^i\widetilde f_*\omega_{\widetilde V/\widetilde W}$ 
is semi-positive by Theorem \ref{semi}. 
By the above inclusion, we can check easily that 
$R^if_*\omega _{V/W}$ is semi-positive. 
We note that $W$ is a curve.  
\endproof

The following theorem is a slight generalization of Theorem \ref{semi}. 
It seems to be well-known to specialists. 
It is buried in Kawamata's proof of the semi-positivity theorem (see 
\cite[\S 4]{ka0}). 
We note that $W$ is not necessarily complete in Theorem \ref{iide} 
below. 

\begin{thm}\label{iide}
Let $f:V\to W$ be a proper surjective morphism 
between non-singular varieties with connected fibers. 
Assume that there exists a simple normal crossing divisor 
$\Sigma$ on $W$ such that $f$ is smooth over $W_0:=W\setminus 
\Sigma$. 
We put $V_0:=f^{-1}(W_0)$, $f_0:=f|_{V_0}$, 
and $k:=\dim f$. 
We assume that all the local monodromies on $R^{k+i}f_{0*}\mathbb C
_{V_0}$ around $\Sigma$ are unipotent. 
Let $C$ be a complete curve on $W$. 
Then the restriction $(R^{i}f_{*}\omega_{V/W})|_{C}$ is semi-positive. 
\end{thm}
\proof
By Chow's lemma, desingularization theorem, and 
\cite[Proposition 1]{ka2}, we can assume that 
$V$ and $W$ are quasi-projective and $\Sigma$ is a simple 
normal crossing divisor. 
We note that if $\psi:W'\to W$ is a proper birational 
morphism from a quasi-projective variety, 
then there exists a complete curve $C'$ on $W'$ 
such that $\psi(C')=C$. 
Let $g:\widetilde C \to C\subset W$ be the normalization of $C$. 
If $C\cap W_0\ne \emptyset$, 
then $g^*R^if_*\omega _{X/Y}\simeq R^ih_*\omega _{D/\widetilde C}$ 
is semi-positive by Theorem \ref{semi}, where $D$ is a desingularization 
of the main component of $V\times _W\widetilde C$ and 
$h:D\to \widetilde C$. 
So, we have to treat the case when $C\subset \Sigma$. 
By cutting $W$, we can take an irreducible surface $S$ on $W$ 
such that $S\not\subset \Sigma$ and $C\subset S$. 
Take a desingularization $\pi:\widetilde S\to S\subset W$ of $S$ 
such that $\pi^{-1}(\Sigma)$ is simple normal crossing on $\widetilde 
S$ and there is a smooth projective 
irreducible curve $\widetilde C$ on $\widetilde
S$ with $\pi(\widetilde C)=C$. 
By \cite[Proposition 1]{ka2}, 
we have $\pi^*R^if_*\omega_{V/W}\simeq 
R^ih_*\omega_{T/\widetilde S}$, where $T$ is a desingularization 
of the main component of $V\times _W\widetilde S$ and 
$h:T\to \widetilde S$. 
So, it is sufficient to check that $\pi^*R^if_*\omega_{V/W}\simeq 
R^ih_*\omega_{T/\widetilde S}$ is 
semi-positive on $\widetilde C$. 
We compactify $h:T\to \widetilde S$. 
Then we obtain an algebraic fiber space $\overline h:\overline T\to 
\overline S$. 
After modifying $\overline h$ birationally, 
we can assume that $\Delta:=(\overline S\setminus \widetilde S)\cup 
\pi^{-1}(\Sigma)$ is a simple normal crossing divisor 
on $\overline S$. 
We can assume that $\overline h$ is smooth over $\overline S\setminus
\Delta$. 
We note that $\widetilde C$ is an irreducible component of $\Delta$. 
By Kawamata's covering trick \cite[Theorem 17]{ka0}, 
we can take a finite cover $\varphi:S'\to \overline S$ which 
induces a unipotent reduction $h':T'\to S'$ 
(see \cite[Corollary 18]{ka0} or \cite[(4.5)]{M}). 
We put $U:=\varphi^{-1}(\widetilde S)$. By Theorem \ref{semi}, 
$R^ih'_*\omega_{T'/S'}$ is semi-positive and 
$(R^ih'_*\omega_{T'/S'})|_{U}\simeq 
\varphi^*
((R^i\overline {h}_*\omega_{\overline {T}/\overline{S}})|_{\widetilde S}) 
\simeq \varphi^*
(R^i{h}_*\omega_{{T}/\widetilde{S}})$. 
Therefore, it is easy to check that 
$(R^i{h}_*\omega_{{T}/\widetilde{S}})|_{\widetilde C}$ is 
semi-positive since $\widetilde C\subset \widetilde S$. 
So, we obtain the required result. 
\endproof

The following corollary is obvious by Theorem \ref{iide}. 
It will be useful when we study algebraic fiber spaces 
in the relative setting.  

\begin{cor}\label{rel}
In the same notation and assumptions as in {\em{Theorem \ref{iide}}}, 
we further assume that $W$ is proper over a variety $B$. 
Then $R^if_*\omega_{V/W}$ is semi-positive over $B$, 
that is, $(R^if_*\omega_{V/W})|_{F}$ is semi-positive 
for every fiber $F$ of $W\to B$. 
\end{cor}

\begin{rem}\label{last}
In the proof of Theorem \ref{iide}, 
we only use the following fact that 
$R^if_*\omega_{V/W}$ is characterized as the canonical 
extension of the bottom Hodge filtration. 
\end{rem}

The next theorem is an easy consequence of the semi-positivity 
theorem. For the proof, we use Theorem \ref{curve}.  

\begin{thm}\label{irr}
Let $f:X\to Y$ be an algebraic fiber space. Then 
$$q(Y)\leq q(X)\leq q(Y)+q(F), 
$$ 
where $F$ is a general fiber. 
Moreover, if $q(F)=0$, then $\Alb(X)\to \Alb(Y)$ is an isomorphism. 
\end{thm}
\proof 
By modifying $f$ birationally, we can assume that 
there exists a simple normal crossing divisor $\Sigma$ on $Y$ 
such that $f$ is smooth over $Y\setminus \Sigma$ and 
$\Supp f^*\Sigma$ is simple normal crossing. 
By Leray's spectral sequence, we have 
$$
0 \to H^1(Y,\mathcal O_Y) \to
H^1(X,\mathcal O_X)\to H^0(Y,R^1f_*\mathcal O_X)\to 0. 
$$ 
We note that $H^2(Y,\mathcal O_Y)\to H^2(X,\mathcal O_X)$ 
is injective. 
Therefore, $q(Y)\leq q(X)=q(Y)+h^0(Y,R^1f_*\mathcal O_X)$. 
Let $H$ be a very ample divisor on $Y$ such that 
$h^0(Y, R^1f_*\mathcal O_X\otimes \mathcal O_Y(-H))=0$. 
We take $H$ general. 
We have; 
$$
0\to R^1f_*\mathcal O_X\otimes \mathcal O_Y(-H)
\to R^1f_*\mathcal O_X \to R^1f_*\mathcal O_{f^{-1}H}\to 0. 
$$ 
Thus, we obtain that $h^{0}(Y,R^1f_*\mathcal O_X)\leq 
h^{0}(Y,R^1f_*\mathcal O_{f^{-1}H})$. 
By repeating this argument, we can assume that 
$Y$ is a curve. 
By Theorem \ref{curve}, 
$(R^i f_*\mathcal O_{X})^{\vee}\simeq 
R^{\dim f-i}\omega_{X/Y}$ 
is semi-positive. 
Then $h^0(Y,R^1f_*\mathcal O_X)\leq q(F)$ 
by Lemma \ref{seminega} below. 
We note that the rank of $R^1f_*\mathcal O_X$ is $q(F)$. 
Thus, we get the required inequality. 
The latter part is obvious. 
\endproof

\begin{lem}\label{seminega}
Let $\mathcal E$ be a vector bundle on a non-singular 
projective curve $V$. 
Assume that the dual vector bundle $\mathcal E^{\vee}$ 
is semi-positive. 
Then $h^0(V,\mathcal E)\leq r$. 
In particular, if $h^0(V,\mathcal E)=r$, then $\mathcal E$ is 
trivial. 
\end{lem}
\proof 
We take a basis $\{\varphi_1,\cdots, \varphi_l\}$ of 
$H^0(V,\mathcal E)$. 
We have to prove that $l\leq k$. 
We define $\psi_i:=\varphi_1\oplus\cdots\oplus\varphi_i:
\mathcal O_V^{\oplus i}\to \mathcal E$ for every $i$. 
We put $\psi_0=0$ for inductive arguments. 

\begin{claim}
$\xIm \psi_i\simeq \mathcal O_V^{\oplus i}$ is a subbundle 
of $\mathcal E$ for every $i$, where $\xIm \psi_i$ denotes 
the image of $\psi_i$. 
\end{claim}
\proof[Proof of Claim]
We use the induction on $i$. 
We assume that $\xIm \psi_{i-1}\simeq \mathcal O_V^{\oplus (i-1)}
$ is a subbundle 
of $\mathcal E$. 
Thus, we have $h^0(V, \xIm \psi_{i-1})=i-1$. 
Therefore, $\xIm \psi_i$ has rank $i$. 
Let $\mathcal F$ be the double dual of $\xIm \psi _i$. 
Then $\mathcal F$ is a rank $i$ subbundle of $\mathcal E$ such 
that $\psi_i$ factors through $\mathcal F$. 
Since $\mathcal F^{\vee}$ is semi-positive by the 
semi-positivity of $\mathcal E^{\vee}$ and 
$\mathcal F$ is semi-positive by the definition of $\mathcal F$, 
we obtain that 
$\psi_i:\mathcal O_V^{\oplus i}\simeq \mathcal F$. 
\endproof
Thus we obtain that $l\leq r$, that is, 
$h^0(V,\mathcal E)\leq r$. The latter part is obvious. 
\endproof

The following corollary is a slight generalization of 
\cite[Corollary 2]{ka0}. 
The proof is obvious by Theorem \ref{irr}.  

\begin{cor}\label{cha}
Let $f:X\to Y$ be an algebraic fiber space. 
Assume that sufficiently general fibers have zero 
Kodaira dimension. 
Then $0\leq q(X)-q(Y)\leq q(F)\leq \dim f$, 
where $F$ is a general fiber of $f$. 
Furthermore, if $q(X)-q(Y)=\dim f$, then general fibers 
are birationally equivalent to Abelian varieties. 
\end{cor}

\ifx\undefined\bysame
\newcommand{\bysame}{\leavevmode\hbox to3em{\hrulefill}\,}
\fi

\end{document}